\documentclass[twocolumn]{article}

\usepackage{arxiv}

\usepackage[utf8]{inputenc} 
\usepackage[T1]{fontenc}    
\usepackage{hyperref}       
\usepackage{url}            
\usepackage{booktabs}       
\usepackage{amsfonts}       
\usepackage{amsmath}
\usepackage{nicefrac}       
\usepackage{microtype}      
\usepackage{cleveref}       
\usepackage{lipsum}         
\usepackage{graphicx}
\usepackage{natbib}
\usepackage{doi}
\usepackage{capt-of}

\usepackage{caption}
\usepackage{subcaption}
\usepackage{svg}
\usepackage{acro}

\DeclareAcronym{PCA}{
  short = PCA,
  long  = Principal Component Analysis
}
\DeclareAcronym{CROCKER}{
  short = CROCKER,
  long  = Contour Realization Of Computed k-dimensional hole Evolutions in the Rips complex
}

\title{A study of holes: Topological analysis reveals crowd dynamics regimes in a bidirectional corridor scenario}

\date{July 6, 2026}

\author{ \href{https://orcid.org/0009-0005-2122-8877}{\includegraphics[scale=0.06]{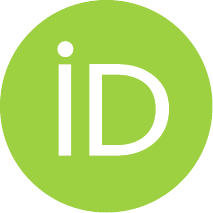}\hspace{1mm}Sabrina Desiree Kern} and \href{https://orcid.org/0000-0002-3369-6206}{\includegraphics[scale=0.06]{figures/orcid.pdf}\hspace{1mm}Gerta K\"oster} \\
	Faculty of Informatics and Mathematics\\
	Hochschule M\"unchen -- University of Applied Sciences\\
	Lothstra{\ss}e 34, 80335 M\"unchen, Germany \\
	\texttt{sabrina.kern@hm.edu} \\
}

\renewcommand{\shorttitle}{A study of holes: Topological analysis reveals crowd dynamics regimes in corridor}

\hypersetup{
pdftitle={A study of holes: Topological analysis reveals crowd dynamics regimes in a bidirectional corridor scenario},
pdfauthor={Sabrina Kern, Gerta  K\"oster},
pdfkeywords={crowd dynamics, topological data analysis, persistent homology, corridor, bidirectional flow
},
}

\begin{document}
\twocolumn[{
  \begin{@twocolumnfalse}
	\maketitle
	\begin{abstract}
		This study harnesses topological analysis in an attempt to reveal structure in the dynamics of a crowd. 
	 Topology and in particular persistent homology characterizes relational structures in data through the number of \textit{connected components} and \textit{holes}, that is, a loop of pairwise connection with no connections across it. We apply this universal data analysis method to a simulated time series of individual pedestrian positions of a crowd moving through a wide corridor -- either uni- or bidirectional. We consider two pedestrians to be connected, when they are sufficiently close. This approach leads to two matrices containing the persistence signatures for the whole time series, so-called \textit{CROCKERs}. 
	 Despite the high level of data abstraction, the CROCKERs' first two principal components on time-delayed positional data show a clear separation of the different parameter configurations. This holds up to symmetry. Our results support our claim that persistent homology is a useful tool to characterize crowd dynamics without introducing any prior assumptions about the detectable spatio-temporal patterns. 
	\end{abstract}

  \end{@twocolumnfalse}
  \vspace{1em}
}]


\section{Introduction: The topological view of pedestrian dynamics}
When looking at a crowd from above, we see areas of closeness and gaps that seem to depend on the topography, but also on individual and collective destinations, and even the psychological state of the crowd. 
Relative positions in a crowd of fans moving quickly towards a concert of their favourite artist will certainly differ from those in a crowd leisurely waiting for a train.
The occurrence of topological structures, that is, connectivity and holes, inspires us to investigate crowd dynamics through the lens of topology.
Topology and in particular \textit{persistent homology} offers a tool to analyze these geometrical structures without exerting any expert knowledge. 
We are encouraged by previous work in \cite{bhaskar-2019-cdyn}, who successfully classify emergent collective phenomena in interaction models from biology. 
Interestingly, their approach does not need any hand-crafted metrics to distinguish different collective patters and to solve the inverse problem of uncovering the parameters they used for their simulated time series data. 
In contrast to works such as \cite{herranz-2026-cdyn}, we plan to use the pedestrians' positions directly, without introducing additional macroscopic observables such as the velocity of the crowd. 

Our aim of this work is to understand when and how collective phenomena in crowds can be identified as structurally different, that is, when and how they can be classified without using hand-crafted metrics. This requires us to transfer established methods from topological data analysis to this new area of application, namely crowd dynamics. 
As we are not aware of other studies using persistent homology to analyze pedestrian crowds through their positions, this work starts on a simple but dynamically interesting example: a corridor with uni- or bidirectional flow. 
This minimal example already contains qualitative phenomena such as lane formation while keeping the number of parameters rather small. 
To better understand how much information is captured in the case of a dynamic crowd, we want to compare several scenarios in a low-dimensional state space, where we expect scenarios with similar spatio-temporal dynamics to cluster. This in turn provides some insight into the capabilities of the topological data analysis in contrast to using hand-crafted metrics.  

\section{Methods: Obtaining macroscopic observables using persistent homology}
\begin{figure}[htp!]
  \centering
  \includegraphics[width=\linewidth, trim={2.2cm 4.5cm 8.3cm 2.2cm}, clip]{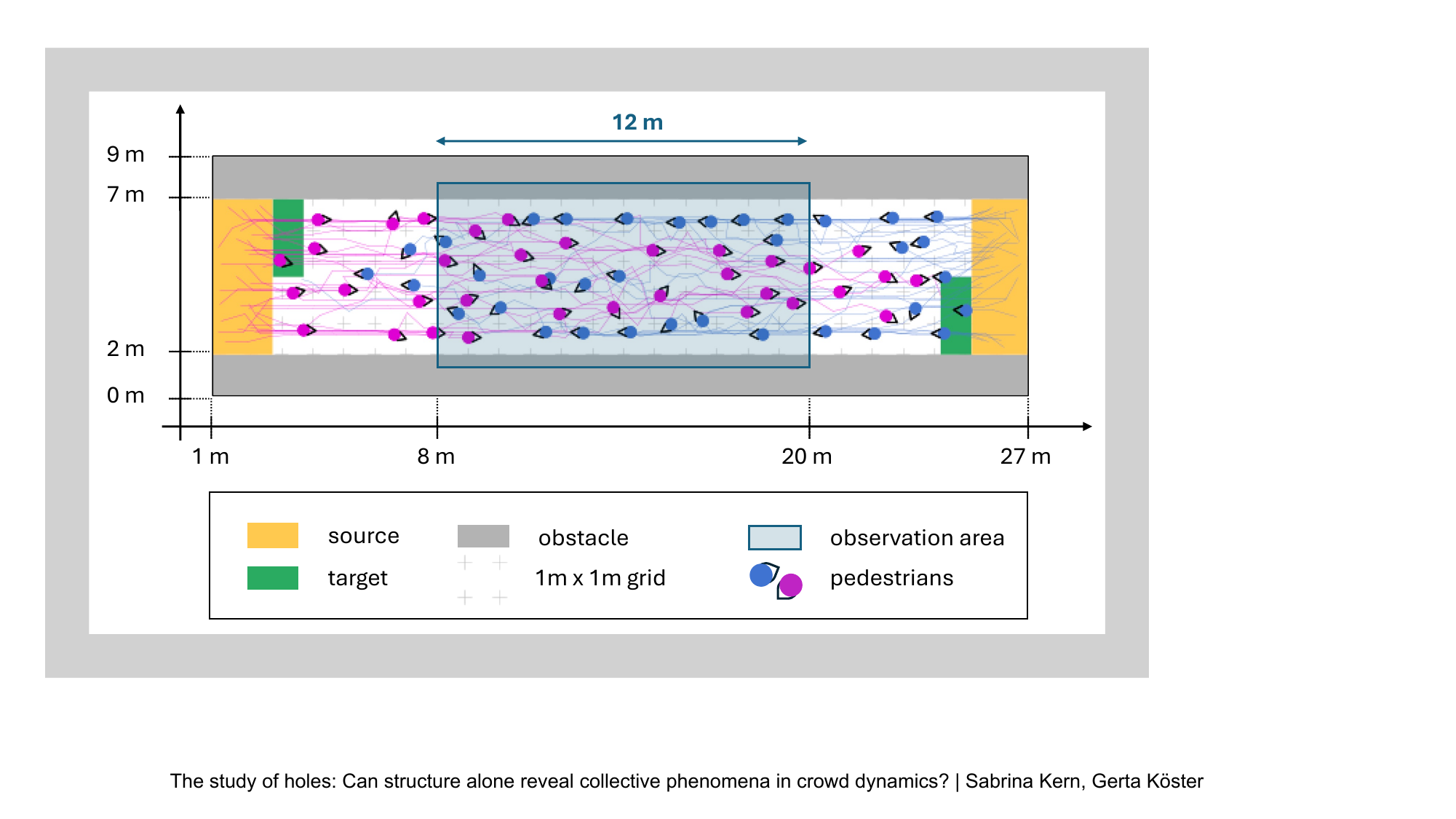}
    \caption{Geometry of the simulated bidirectional corridor scenario: Source and targets are used on both ends, with a width of 5 m and an observation area of 12 m length. Pedestrians appear out of the outer source areas and disappear once reaching the target areas.}
    \label{fig:vadere-layout}
\end{figure}

The aim of this study is to showcase which phenomena topological data analysis is able to capture and which not. 
We create our data set using the open-source crowd simulator Vadere \cite{kleinmeier-2019-cdyn} to 
allow full control over the boundary conditions, while the corridor geometry is inspired by the laboratory experiments described in \cite{zhang-2012b-cdyn}. We model the corridor at a length of 26 meters and a width of 5 meters, as shown in Figure \ref{fig:vadere-layout}. The observation area for which we extract the pedestrian positions is placed in the middle of the corridor with a width of 6 meters and length of 12 meters. This ensures that the data set only includes the locomotion dynamics without any artifacts due to appearing and disappearing of pedestrians at the sources and targets. 
The locomotion is modelled by the \textit{Optimal Steps Model} for which we see the occurrence of lanes, highlighted by colouring based on the given target (left or right). The forming of lanes for pedestrians walking into the same direction is fostered by the shift of the targets to the right hand side of each crowd. 

For our corridor scenario, we vary the total inflow between 4, 6, an 8 pedestrians every two seconds, divided between the two sources. This induces three different regimes: the unidirectional flow (all persons from one side, none from the other), balanced bidirectional flow (half of the persons from each side) and all combinations of unbalanced flows ($n$ persons from one side and $m-n$ from the other side, $m = 4, 6, 8$). 
We always include the symmetric configuration, resulting in a total of 21 different combinations to which we will refer to as \textit{scenarios}. 
We label the scenarios by left and right source inflow, for example \textit{L01R07} describing the scenario with 1 person every 2 seconds created on the \underline{L}eft hand side and 7 on the \underline{R}ight hand side respectively. For each of the parameter configurations, we simulate for 50 different initializations of locomotion parameter choices in the simulation which are considered stochastic; for example the preferred walking speed for each simulated pedestrian is drawn from a truncated normal distribution.
This leads to slightly different trajectories for each simulation run. 

\begin{figure}
    \centering
    \includegraphics[width=1\linewidth, trim={6.7cm 1.8cm 15.5cm 2.2cm}, clip]{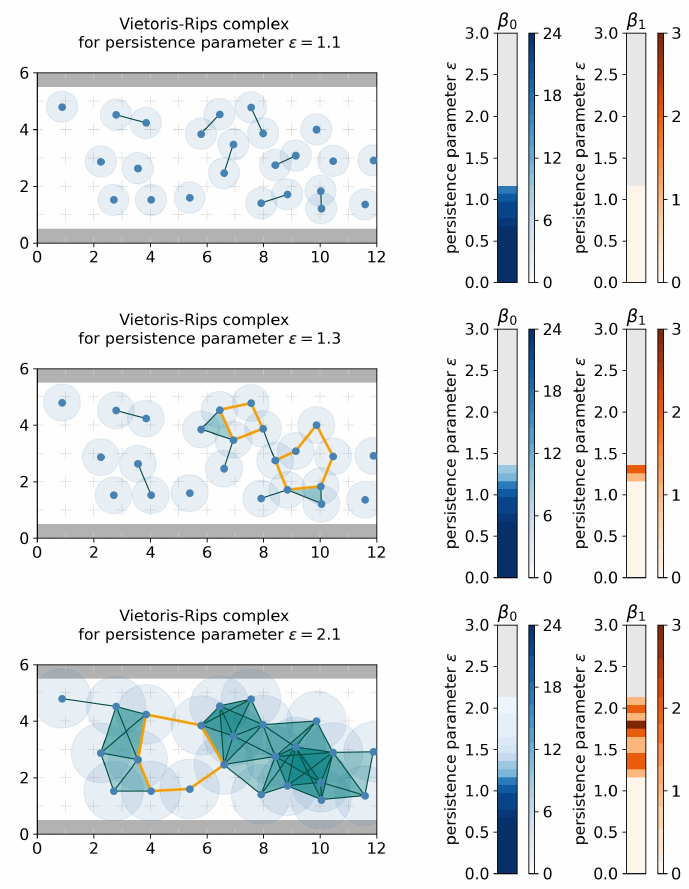}
    \caption{Vietoris-Rips complex for three choices of the persistence parameter $\epsilon$ on the same positional arrangement of the crowd. The two bars on the right show the corresponding Betti numbers $\beta_0$ and $\beta_1$ in dependency with the change of $\epsilon$. }
    \label{fig:persistence-parameter-change}
\end{figure}

Following the promising results in \cite{bhaskar-2019-cdyn}, we proceeded as follows for each individual simulation:
We represent pedestrians by their positions as points $(x_t, y_t)$ on a 2-dimensional plane, with one snapshot at time $t$ every 0.2 seconds. 
Given a fixed threshold -- the \textit{persistence} or \textit{scale parameter} $\epsilon$ --, we construct a graph-like structure via connecting those pedestrians which are close to one another. Here, we use the Euclidean distance to measure the closeness, introducing an edge when their distance is smaller than $\epsilon$ and fill out triangles for each set of three mutually-close points. 
The emerging object is referred to as the \textit{Vietoris–Rips complex}, displayed in Figure \ref{fig:persistence-parameter-change} for different persistence parameters $\epsilon$. The Vietoris–Rips complex contains the 0-, 1-, and 2-simplices corresponding to the pedestrians (nodes), their connections (edges) and the triples with pairwise connections respectively (triangles). Put into more general terms, a \textit{k-simplex} is defined as the complex hull of $(k+1)$ affinely independent points. 
Any cycle in the Vietoris-Rips complex between more than three nodes, which does not contain any 1-simplex, is referred to as a 1-dimensional \textit{hole}. Examples of holes are highlighted in the two lower panels in Figure \ref{fig:persistence-parameter-change} for $\epsilon=1.5$. 

At each value of $\epsilon$ -- the threshold of being close --, the topological structure among pedestrian positions is then characterized by its number of connected components, \textit{Betti number $\beta_0$}, and 1-dimensional holes, \textit{Betti number $\beta_1$}, in the Vietoris–Rips complex. 
If one varies the persistence parameter $\epsilon$ the corresponding Betti numbers change as well. Variation in discrete steps leads to one vector per Betti number and can be understood as a signature for a specific moment in time, thus capturing the persistence (emergence and disappearance) of connected components and distinct holes. 
These signatures are depicted on the right hand side of Figure \ref{fig:persistence-parameter-change}; clearly showing that the number of connected components $\beta_0$ monotonically decreases with $\epsilon$ while the number of holes for this specific set of points peaks at $\epsilon=1.8$ with 3 holes. For the visualization, we discretize $\epsilon$ in steps of 0.1, showing 3 different $\epsilon$ values alongside their $\beta_0$ and $\beta_1$ vectors up to this point. 
For a more detailed introduction into persistent homology, the reader is referred to \cite{ghrist-2007-math}. 

\begin{figure*}
    \centering
    \includegraphics[width=\linewidth, trim={1.5cm 2.7cm 3.2cm 1cm}, clip]{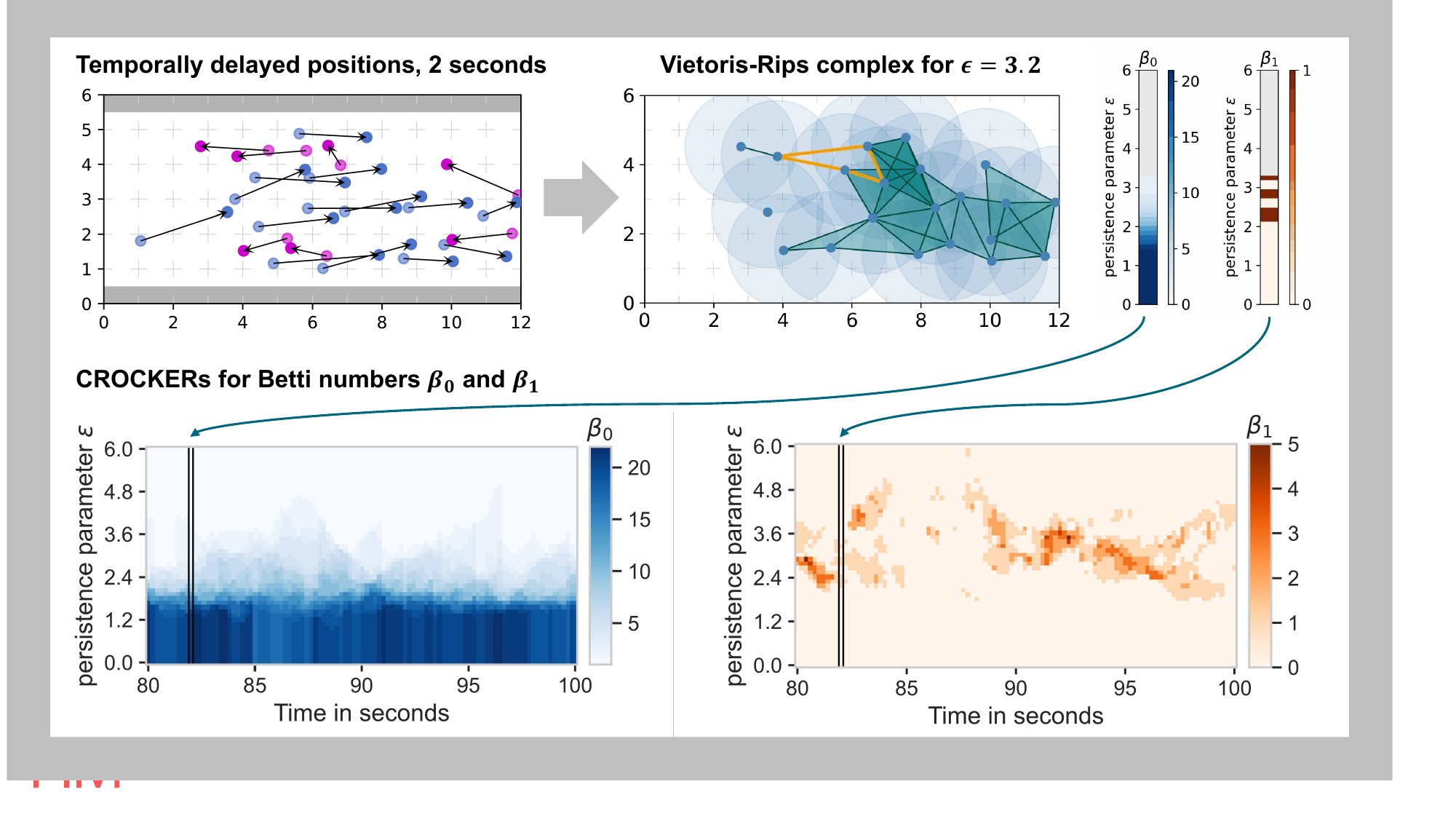}
    \caption{\acp{CROCKER} for a single time series with temporal delay of 2 seconds: In the upper left, the time delay is shown for a single frame, keeping only pedestrian which are in the observation area for the current and the delayed frame. In the upper right, we show the corresponding Vietoris-Rips complex (using only the leading position for visualization of the simplices) and its Betti numbers $\beta_0$ and $\beta_1$ up to $\epsilon = 3.2$. The frame's signatures differ significantly from the one shown in Figure \ref{fig:persistence-parameter-change} where no temporal delay is used on the same relative positions.
    The \acp{CROCKER} in the lower panel summarize the values of the Betti numbers for the whole time series, each column representing the persistence of connected components (left) and holes (right) for a single positional snapshot. The rows correspond with different (ordered) choices of epsilon. For each time series we only consider the 20 seconds in the near steady flow interval between 80 and 100 seconds. }
    \label{fig:crockers-for-ts}
\end{figure*}

So far, the positions in each frame are treated individually, without any directional information. 
To also include directional information for each simulated pedestrian, we use a temporal delay embedding: to each pedestrian position we add its predecessor. The data points per pedestrian now take the form $(x_t, y_t, x_{t-d}, y_{t-d})$ with $d$ denoting the delay in time $t$ in seconds. Adding the delay increases the input dimensions from 2 to 4, inducing a different geometry in the positional information of each time step and therefore leading to different signatures of the Betti numbers as well (compare upper two panels in Figure \ref{fig:crockers-for-ts}). 
 
Given the persistence information for all snapshots in each simulated time series  -- either for the 2-dimensional positions or the time-delayed ones --, we can summarize them into a single object, the so-called \textit{\ac{CROCKER} plot},
with one column per snapshot and one row per discrete choice of $\epsilon$ \cite{topaz-2015-math}. 
Both of the \acp{CROCKER} for one exemplary time series are shown in the lower part of Figure \ref{fig:crockers-for-ts}.
Now, these two \acp{CROCKER} -- put in a single long vector -- provide a highly abstracted representation of the structural properties in the time series of pedestrians moving through a corridor.

To compress the information contained in the \acp{CROCKER}, we apply \ac{PCA} such that each simulated time series is now represented by a single data point in a low-dimensional space. This approach allows us to directly observe potential clusters for different boundary conditions. In the work by \cite{bhaskar-2019-cdyn}, different regimes clearly separated in the first three principal components. 
In our study, we quantify the separation for the scenarios as clusters using \texttt{sklearn}'s silhouette coefficient as a score \cite{rousseeuw-1987-math}. The clearer the separation and the denser the obtained clusters, the higher the score which always ranges between -1 and 1. Furthermore, we compare the resulting scatter diagram visually. 

We use the open-source python library \texttt{ripser} \cite{ctralie-2018-cs} to create the \acp{CROCKER}, as well as \texttt{scipy} and \texttt{sklearn} for statistical analysis. The corresponding code 
can be shared upon request. 

\section{Results: Persistent homology distinguishes simulated counterflow regimes}
\begin{figure*}
     \centering
     \begin{subfigure}[b]{0.35\textwidth}
         \centering
         \includegraphics[height=6.4cm, trim={2cm 3cm 20cm 3.5cm}, clip]{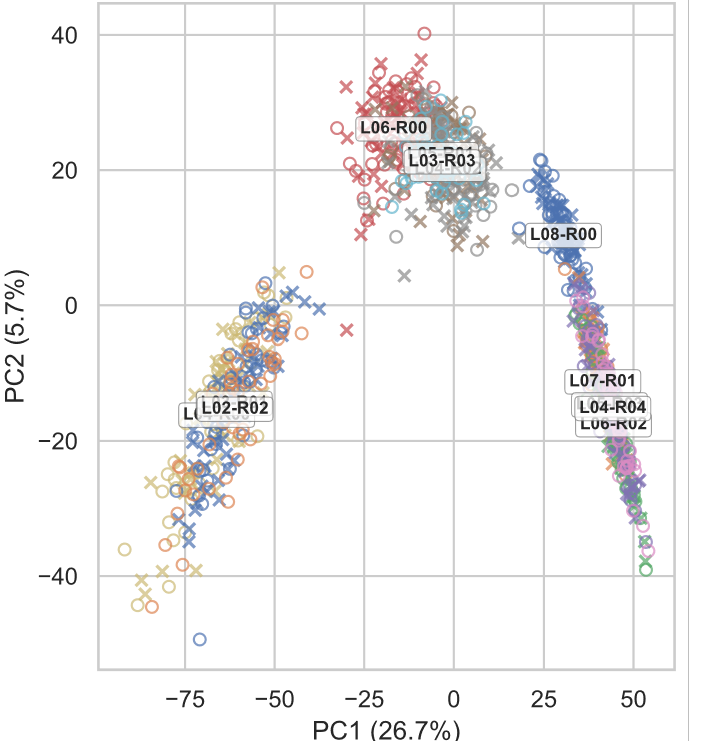}
         \caption{\acp{CROCKER} based on 2d-positions.}
         \label{fig:pca-no-delay}
     \end{subfigure}
     \hfill
     \begin{subfigure}[b]{0.6\textwidth}
         \centering
         \includegraphics[height=6.4cm, trim={14cm 3cm 2.3cm 3.5cm}, clip]{figures/pca_time_delay.pdf}
         \caption{\acp{CROCKER} based on 4d-positions with a temporal delay of 2 seconds.}
         \label{fig:pca-with-delay}
     \end{subfigure}
        \caption{\ac{PCA} for the simulated time series data with one data point representing one time series: We fit two \acp{PCA} on the whole data set, one using a temporal delay (right) and one without (left). In both clusterings the different scenarios can be separated, while our approach of persistent homology does not differentiate the symmetric scenarios (when swapping \underline{L}eft and \underline{R}ight for the defined inflow). Using the temporal delay leads to a significantly better separation of the different regimes, with a silhouette coefficient of 0.376 instead of 0.033.}
        \label{fig:pca-clustering}
\end{figure*}

For each series, and each snapshot, we perform a topological analysis as described above, that is, we vary the persistence parameter $\epsilon$ in 50 steps between 0 and 6 (meters). We summarize persistence properties for the connected components ($\beta_0$) and 1d-holes ($\beta_1$) into \acp{CROCKER}; one with and one without temporally delayed positions.
Once we have obtained our two sets of \ac{CROCKER}-based representations, we fit one \ac{PCA} each, eventually arriving at two low dimensional representations of all simulations. Now, one data point represents one time series in a two-dimensional space; as visualized in Figure \ref{fig:pca-clustering}. 

Figure \ref{fig:pca-clustering} also compares the effect of using the temporal delay before constructing the \acp{CROCKER}: 
the representation in the left panel is obtained without a temporal delay, while the representation on the right is based on a temporal delay of 2 seconds. 
In both cases, we clearly see clusters. In Panel \ref{fig:pca-no-delay}, one cluster forms for each total inflow. However, we see no separation that reveals how many agents enter the corridor from each side. 
In contrast to this, in Panel \ref{fig:pca-with-delay}, more clusters appear that do not only correspond to the total inflow but also 
indicate the different regimes of uni- and bidirectional flow. 
The clusters emerge up to symmetry (circle and x-markers), separating scenarios of (un-)balanced inflows but not considering which side of the corridor has more entries (left or right). Indeed, intuitively symmetry only changes the point of observation for the same dynamics in the relative positions and shows that our abstract method is invariant to rotation. This is a significant advantage over grid-based approaches of macroscopic observables such as the density; in particular when structure and dynamic is of main interest, independent of the exact topography.  

In Figure \ref{fig:pca-clustering}, the temporal delay of $d = 2$ seconds, leads to a silhouette coefficient of 0.376 for the symmetric scenarios in the first two principal components. This is the best score among the choices of $d = 0.8, 1.2, 1.6, 2, 2.4$ seconds. Without a temporal delay the silhouette coefficient is as low as 0.033.
Intuitively, the temporal delay implicitly includes dynamical information for each data point, providing the direction of movement in addition to the position. Overall, Figure \ref{fig:pca-clustering} clearly shows the importance of the temporal delay in our analysis. 

By construction, there is a significant correlation between the value of the Betti-numbers -- and as such between the \acp{CROCKER} -- and the number of points or pedestrians in our case. 
In particular, for a very small persistence parameter $\epsilon$, the Betti number $\beta_0$ counting the number of connected components is equal to the number of pedestrians in the observation area. 
Furthermore, the higher the average density the more edges are added with a relatively small persistence parameter value, reducing the number of connected components and hence the Betti number $\beta_0$ early on. 
This being said, given a fitted \ac{PCA} for the \acp{CROCKER}, we check for the correlation between the first two principal components and the average pedestrian count of the corresponding scenario. The first principal component is highly correlated with the pedestrian count with a Pearson coefficient of $\rho_{1} \approx 0.92$, while the second principal component has a weak correlation with $\rho_{2} \approx 0.24$. 
%
%
We are encouraged by the observation that \acp{CROCKER} -- while being a purely structural representation -- confirm the known fact 
that the local density has a significant impact on crowd dynamics, both quantitatively and qualitatively. In future work we want to better understand which information we gain by using structural methods over classic analysis. 

\section{Conclusion}

In conclusion, topological analysis enabled a highly accurate classification of simulation scenarios in a corridor with or without counterflow. This first result encourages us to extent the analysis to more complex topographies and to apply it to data from laboratory experiments or field observations. Furthermore, we hope to inspect the correlation and interaction with density more closely in future works.

\section{Acknowledgements}
This work was financially supported by Deutsche Forschungsgemeinschaft (DFG, German Research Foundation) through the project Idefiks -- Identifikation von Elementarszenarien der Dynamik von Fußgängerströmen über interpretierbare Künstliche Intelligenz (project no. 524081074). 

\bibliographystyle{plain}
\bibliography{literature}

\end{document}